\documentclass[12pt]{article}
\usepackage{amsmath,mathrsfs}
\usepackage{amssymb}
\usepackage{latexsym}
\usepackage{amsfonts}

\newcommand{\qed}{\ifmmode$\Box$\else{\unskip\nobreak\hfil
\penalty50\hskip1em\null\nobreak\hfil$\Box$
\parfillskip=0pt\finalhyphendemerits=0\endgraf}\fi}

\newcommand{\ra}{\rangle}
\newcommand{\la}{\langle}

\newcommand{\cF}{{\mathcal F}}

\newcommand{\cH}{\mathcal H}

\newcommand{\C}{\mathbb C}
\newcommand{\D}{\mathbb D}
\newcommand{\T}{\mathbb T}
\newcommand{\z}{{\mathfrak z}}
\newcommand{\w}{\omega}
\newtheorem{Pa}{Paper}[section]
\newtheorem{Tm}[Pa]{{\bf Theorem}}
\newtheorem{La}[Pa]{{\bf Lemma}}

\newtheorem{Pn}[Pa]{{\bf Proposition}}

\newtheorem{Dn}[Pa]{{\bf Definition}}

\newcommand{\CC}{{\mathchoice
{\setbox0=\hbox{$\displaystyle\rm C$}\hbox{\hbox
to0pt{\kern0.4\wd0\vrule
height0.9\ht0\hss}\box0}}
{\setbox0=\hbox{$\textstyle\rm C$}\hbox{\hbox
to0pt{\kern0.4\wd0\vrule
height0.9\ht0\hss}\box0}}
{\setbox0=\hbox{$\scriptstyle\rm C$}\hbox{\hbox
to0pt{\kern0.4\wd0\vrule
height0.9\ht0\hss}\box0}}
{\setbox0=\hbox{$\scriptscriptstyle\rm
C$}\hbox{\hbox to0pt{\kern0.4\wd0\vrule
height0.9\ht0\hss}\box0}}}}

\begin{document}
\title{A characterization of Schur multipliers between
character-automorphic Hardy spaces }
\author{D. Alpay and M. Mboup}
\date{}
\maketitle
\begin{abstract}
We give a new characterization of
character-automorphic Hardy spaces of order $2$
and of their contractive multipliers in terms of
de Branges Rovnyak spaces. Keys tools in our
arguments are analytic extension and a
factorization result for matrix-valued analytic
functions due to Leech.
\end{abstract}
{\bf Keywords.} Character-automorphiuc functions,
Hardy spaces, de Branges Rovnyak spaces, Schur
multipliers.\\

{\bf Mathematics Subject Classification (2000).}
Primary: 30F35, 46E22. Secondary: 30B40
\section{Introduction}
\setcounter{equation}{0} Let $\Gamma$ be a
Fuchsian group of M{\"o}bius transformations of
the unit disk $\D = \{z \in \C\, ; \, |z| < 1\}$
onto itself. For $1 \leqslant p \leqslant
\infty$ and for any character $\alpha$ of
$\Gamma$,  we consider the spaces
$${\mathcal H}_p^\alpha =
\left\{f \in {\mathcal H}_p\, \, \mid \, f\circ
\gamma = \alpha(\gamma) f, \quad \forall \,
\gamma \in \Gamma\right\}.$$
These spaces are called character-automorphic
Hardy spaces. A characterization of such spaces
in terms of Poincar\'e theta series may be found
in \cite{Pommerenke}, \cite{Samokhin},
\cite{metzger}, \cite{Marden}. In particular,
Pommerenke showed in  \cite{Pommerenke} that the
series
\begin{equation}
\label{poincare}
  f(z) =
  \frac{b_0(z)}{b'_0(z)}\sum_{\gamma \in \Gamma}
  \overline{\alpha(\gamma)} \theta(\gamma(z))
  h(\gamma(z)) \frac{\gamma'(z)}{\gamma(z)}
\end{equation}
defines a bounded linear operator from the
classical Hardy space ${\mathcal H}_p(\D)$ into
the
subspace ${\mathcal H}_p^\alpha(\D)$.\\

In the present paper we restrict ourselves to
the case $p=2$. We first give a characterization
of the character-automorphic Hardy space
${\mathcal H}_2^\alpha({\mathbb D})$ in terms of
an associated de Branges Rovnyak space of
functions analytic in the open unit disk; see
Theorem \ref{hardy}. We also characterize the
contractive multipliers between ${\mathcal H
}_2^{\overline{\beta}\alpha}({\mathbb D})$ and
${\mathcal H }_2^\alpha({\mathbb D})$, where
$\alpha$ and $\beta$ two given characters; see
Theorem \ref{schur}. Our method is mainly based
on analytic extension of positive kernels and
factorization results from Nevanlinna-Pick
interpolation theory.

\section{A review on character-automorphic
Hardy spaces}
\subsection{Fuchsian groups and automorphic functions}
Let $G$ be a group of linear transformations,
$T(z) = \frac{a z + b}{c z
 + d},\,\, a d - b c = 1$, in the complex plane
and let $\iota$ denotes the identity
transformation. Two points $z$ and $z'$ in $\C$
are said to be \textit{congruent} with respect
to $G$, if $z'=T(z)$ for some $T \in G$ and $T
\neq \iota$. Two regions $R, R' \subset \C$ are
said to be $G$-congruent or $G$-equivalent if
there exists a transformation $T \neq \iota$
which sends $R$ to $R'$. A region $R$ which does
not contain any two $G$-congruent points and such
that the neighborhood of any point on the
boundary contains $G$-congruent points of $R$ is
called a \textit{fundamental region} for $G$. A
\textit{properly discontinuous} group  is a
group $G$ having a fundamental region
\cite{Ford}. This amounts to saying that the
identity transformation is isolated.
\begin{Dn}
  A \textit{Fuchsian group} is a properly discontinuous group each of whose
  transformation maps $\D$, $\T$ and $\C\backslash\overline{\D}$ onto
  themselves.
\end{Dn}
A Fuchsian group $\Gamma$ is said to be of
convergence type (see \textit{e.g.}
\cite{Pommerenke}) if
$$\sum_{\gamma \in \Gamma} \left(1 - |\gamma(z)|^2\right) = \left(1 -
 |z|^2\right)\sum_{\gamma \in \Gamma} |\gamma'(z)| < \infty \quad z \in \D.
$$
Then, the Green's function \cite{Pommerenke} of
$\Gamma$ with respect to a point $\xi \in \D$ is
defined as the Blaschke product
\begin{equation}\label{Green}
 b_{\xi}(z) = \prod_{\gamma \in \Gamma}
 \frac{\gamma(\xi) - z}{1 - \overline{\gamma(\xi)} z} \frac{|\gamma(\xi)|}{\gamma(\xi)}.
\end{equation}
It satisfies
\begin{equation}\label{charaut}
b_{\xi}(\varphi(z)) = \mu_{\xi}(\varphi)
b_{\xi}(z), \quad \forall \varphi \in \Gamma,
\end{equation}
where $\mu_{\xi}$ is the character of $\Gamma$
associated with $b_{\xi}(z)$. A function
satisfying the relation \eqref{charaut} is said
to be \textit{character-automorphic} with
respect to $\Gamma$ while a $\Gamma$-periodic
function, as for example $|b_{\xi}(z)|=
|b_{\xi}(\varphi(z))|$, is called
\textit{automorphic} with respect to $\Gamma$.
\subsection{Spaces of character-automorphic functions}
We now briefly mention the main properties
pertaining to spaces of character-automorphic
functions. The materials presented here are
essentially borrowed from \cite{Pommerenke} and
\cite{SodYudit97} (see also \cite{Hasumi,
Samokhin}). Let $\widehat{\Gamma}$ be the dual
group of $\Gamma$, \textit{i.e.} the group of
(unimodular) characters. For an arbitrary
character $\alpha \in \widehat{\Gamma}$,
associate the subspaces of the classical space
$L_2(\T)$
\begin{align*}
L_2^{\alpha} &= \{f \in L_2\,|\,f\circ\gamma =
\alpha(\gamma) f, \,
\forall \gamma \in \Gamma\}\\
\cH_2^\alpha(\D)&=L_2^{\alpha}\bigcap \cH_2(\D)
\end{align*}
Let $\Gamma$ be a Fuchsian group
 without elliptic and parabolic element.
 We say that $\Gamma$ is of
 Widom type if, and only if, the derivative of
 $b_0(z)$ is of bounded characteristic.
 In this case, Widom \cite{Widom} has shown
that the space $\cH^\alpha_\infty$ is not
trivial for any character $\alpha \in
\widehat{\Gamma}$ and we have
\begin{Tm}[Pommerenke \cite{Pommerenke}]
Let $\Gamma$ be of Widom type and let
$\theta(z)$ be the inner factor of $b'_0(z)$. If
$\alpha$ is any character of $\Gamma$ and if
$h(z)$ is in $\cH_p(\D), 1 \leqslant p \leqslant
\infty$, then the function defined by
\eqref{poincare} is in $\cH_p^\alpha(\D)$ and
$$\|f\|_p \leqslant \|h\|_p, \quad f(0) = \theta(0) h(0).$$
\end{Tm}

The Poincar\'e series \cite{Poincare} in
\eqref{poincare} thus defines, in particular, a
projection: $P^\alpha : \theta \cH_2(\D) \to
\cH_2^\alpha(\D)$.  An important property of the
space $\cH_2^\alpha(\D)$ that we will need is
that, point evaluation $f \mapsto f(\xi), \xi
\in \D$ is a bounded linear functional. The
space therefore admits a reproducing kernel
$k^\alpha$:
\begin{align*}
\label{repker} &\la f(z), k^\alpha(z, \xi)
\ra_{\cH_2^\alpha(\D)}
 = f(\xi), \, \xi \in \D
 \quad \text{ for all } f \in
\cH_2^\alpha(\D)\\
&\text{with } k^\alpha(z, \xi) \in
\cH_2^\alpha(\D), \quad \forall \, \xi \in \D
\end{align*}
Since $\cH_2^\alpha(\D) \neq \{const\}$, we have
$k^\alpha(\xi, \xi) = \|k^\alpha(\cdot,
\xi)\|^2_{\cH_2^\alpha(\D)} > 0$ for every $\xi
\in \D$.
In the sequel, the Green's function $b_0(z)$
with respect to $0$,  will be denoted by $b(z)$
for short.

Let $\Gamma$ be a group of Widom type and let
$E$ be associated\footnote{See
\cite{Yuditskii2001} for an example of a
construction of a group $\Gamma$ associated with
a finite union of disjoint arcs of the unit
circle.} to it in such a way that
$\overline{\C}\backslash E$ be equivalent to the
Riemann surface $\D/\Gamma$, obtained by
identifying $\Gamma$-congruent points. Then,
there exists a universal covering map $\z : \D
\to \overline{\C}\backslash E \simeq \D/\Gamma$
such that
\begin{itemize}
  \item $\z$ maps $\D$ conformally onto $\overline{\C}\backslash E$,
\item $\z$ is automorphic with respect to $\Gamma$: $\z \circ \gamma = \z, \quad \forall \, \gamma \in \Gamma$
\item and $\z(z_1) = \z(z_2) \Rightarrow \exists \gamma \in \Gamma \mid z_1 = \gamma(z_2)$
\end{itemize}
In particular, $\z$ maps one-to-one the
\textit{normal fundamental domain} of $\Gamma$
with respect to the origin,
\begin{equation}\label{funddomain}
\cF = \{z \in \D\,:\, |\gamma'(z)| < 1 \quad
\text{for all } \gamma \in  \Gamma, \gamma \neq
\iota\}
\end{equation}
conformally onto some sub-domain of
$\overline{\C}\backslash E$. We assume that $\z$
is normalized so that $(\z b)(0)$  is real and
positive. In all the sequel, the character
associated to the Green's function $b(z)$ will
be denoted by $\mu$. The starting point of the
next section is the following result:
\begin{La}[\cite{Yuditskii97}]
 The reproducing kernel for the space $\cH_2^\alpha(\D)$ has the form
\begin{equation}
\label{noyau} k^\alpha(z, \omega)=c(\alpha)
\frac{
\frac{k^{\alpha\mu}(z, 0)}{b(z)}
k^\alpha(\omega, 0)^*-
\left(\frac{k^{\alpha\mu}(\omega,
0)}{b(\omega)}\right)^*k^\alpha (z,
0)}{\z(z)-\z(\omega)^*}
\end{equation}
where
\begin{equation}
c(\alpha)=\frac{\z(0)b(0)}{k^{\alpha\mu}(0, 0)}
> 0.
\end{equation}
\end{La}
\section{An associated de Branges-Rovnyak space}
\setcounter{equation}{0} In this section we give
a characterization of the space ${\mathcal
H}_2^\alpha({\mathbb D})$ in terms of an
associated de Branges Rovnyak space. To begin,
let
\[
\Omega_+=\left\{z\in{\mathbb D}\,;\, {\rm
Im}~\z(z)>0\right\}.
\]
Setting
\[
\begin{split}
A^\alpha(z)&=\sqrt{\frac{c(\alpha)}{2}}\left(
\frac{k^{\alpha\mu}(z, 0)}{b(z)}+i
k^\alpha(z, 0)\right)\\
B^\alpha(z)&=\sqrt{\frac{c(\alpha)}{2}}\left(
\frac{k^{\alpha\mu}(z, 0)}{b(z)}-i k^\alpha(z,
0)\right)
\end{split}
\]
we can rewrite the reproducing kernel $k^\alpha$
as
\begin{equation}
k^\alpha(z, \omega)=\frac{2A^\alpha(z)}{1-i\z(z)}
\frac{1-S_\alpha(z)S_\alpha(w)^*}{1-\sigma(z)\sigma(w)^*}
\frac{2A^\alpha(w)^*}{1+i\z(w)^*} \label{noyau1}
\end{equation}
where $S_\alpha(z)=B^\alpha(z)/A^\alpha(z)$ and
$\sigma(z)=\frac{1+i\z(z)} {1-i\z(z)}$. We note
that  the functions $A^\alpha(z)$ and
$B^\alpha(z)$ are character-automorphic with the
same character $\alpha$ while $S_\alpha$ and
$\sigma$ are automorphic functions. From now on,
the notation $f^\nu(z)$ will means that the
function $f^\nu(z)$ is character-automorphic
with the superscript $\nu \in \widehat{\Gamma}$
being the associated character, and the notation
$f_\nu(z)$ will stand for a function depending
on the character
$\nu$ (automorphic or not).\\

\begin{Pn}
\label{S0} There exists a Schur function
${\mathscr S}_\alpha$ such that
$S_\alpha(z)={\mathscr S}_\alpha(\sigma(z))$.
\end{Pn}
{\bf Proof:} Since the kernel $k^\alpha(z,
\omega)$ is positive in $\D$, and hence in
$\Omega_+$, it is clear that
$\frac{1-S_\alpha(z)S_\alpha(w)^*}{1-\sigma(z)\sigma(w)^*}$
is also positive in $\Omega_+$. Now, observe
that the function $\sigma$ maps $\Omega_+ \cap
\cF$ into some subset $\Delta \subset \D$ and
this mapping is one-to-one. Let $\varsigma$ be
given by: $(\varsigma \circ \sigma) (z) = z, \,
\forall z \in \Omega_+ \cap \cF$ (in particular
this will also hold for any region congruent to
$\Omega_+\cap \cF$) and let the function
$\widetilde{{\mathscr S}_\alpha}$ be defined on
$\Delta$ by:
$$\widetilde{{\mathscr S}_\alpha}(\lambda) = (S_\alpha\circ \varsigma)(\lambda), \, \forall \lambda \in \Delta.$$
Then it comes that the kernel
$$\frac{1-\widetilde{{\mathscr S}_\alpha}(\lambda)\widetilde{{\mathscr S}_\alpha}(\mu)^*}{1-\lambda \mu^*}$$
is positive on $\Delta$. Now, this implies (see
for instance \cite[Theorem 2.6.5]{MR99g:47016})
the existence of a unique extension of
$\widetilde{{\mathscr S}_\alpha}(\lambda)$,
analytic and contractive in all $\D$. We
subsequently call ${\mathscr S}_\alpha(\lambda)$
this extension, and denote by ${\mathcal
H}({\mathscr S}_\alpha)$ the reproducing kernel
Hilbert space with reproducing kernel
\[
K_{{\mathscr S}_\alpha}(\lambda,\mu)=
\frac{1-{\mathscr S}_\alpha(\lambda){\mathscr
S}_\alpha(\mu)^*} {1-\lambda \mu^*}\]

By construction, the equality
$$S_\alpha(z) = {\mathscr S}_\alpha(\sigma(z))$$
holds for all $z \in \Omega_+ \cap \cF$. Since
$S_\alpha(z)$ is analytic in $\D$, it must also
hold for all $\D$.
\mbox{}\qed\mbox{}\\

In connection with the previous proposition and
the next theorem, we recall that reproducing
kernel Hilbert spaces ${\mathcal H}({\mathscr
S})$ of functions analytic in the open unit disk
and with a reproducing kernel of the form
\[
\frac{1-{\mathscr S}(\lambda){\mathscr
S}(\mu)^*}{1-\lambda\mu^*}
\]
were introduced and studied by de Branges and
Rovnyak; see \cite[Appendix]{dbr1}, \cite{dbr2}.
We also refer the reader to \cite{MR90g:47003}
and \cite{MR99g:47016} for more information on
these and on related spaces.
\begin{Tm}
\label{hardy} The character-automorphic Hardy
space ${\mathcal
   H}^\alpha_2({\mathbb D})$ can be described as
\begin{equation}
\label{hardy1} {\mathcal H}_2^\alpha({\mathbb
D}) =\left\{F(z)=\frac{\sqrt{2}A^\alpha(z)}
{1-i\z(z)}f(\sigma(z))\,\,; f\in{\mathcal
H}({\mathscr S}_\alpha)\right\}
\end{equation}
with the norm
\[
\|F\|_{{\mathcal H}_2^\alpha({\mathbb
 D})}=\|f\|_{{\mathcal H}({\mathscr S}_\alpha)}.
\]
\end{Tm}
{\bf Proof:}  Recall that the map which to $F\in
{\mathcal
 H}_2^\alpha({\mathbb D})$
associates its restriction
 $F|_{\Omega_+}$ to $\Omega_+$ is an isometry from
${\mathcal
 H}_2^\alpha({\mathbb D})$ onto the reproducing kernel Hilbert space
with reproducing $k^\alpha(z,\omega)$ defined by
\eqref{noyau}, where $z,\omega$ are now
restricted to $\Omega_+$. We denote this last
space by ${{\mathcal H}_2^\alpha({\mathbb
D})}\big|_{\Omega_+}$. By Proposition \ref{S0}
and using \eqref{noyau1} we see that the
operator of multiplication by
$\frac{2A^\alpha(z)}{1-i\z(z)}$ is an isometry
from the reproducing kernel Hilbert space
${\mathcal H}$  with reproducing kernel
\[
\frac{1-{\mathscr S}_\alpha(\sigma(z)){\mathscr
S}_\alpha(\sigma(w)^*} {1-\sigma(z)\sigma(w)^*}\]
onto ${{\mathcal H}_2^\alpha({\mathbb
D})}\big|_{\Omega_+}$. Furthermore, the
composition map by $\sigma$ is an isometry from
the de Branges Rovnyak space ${\mathcal
H}({\mathscr S}_\alpha)$
onto ${\mathcal H}$. We have that
\[
{\mathcal H}=\left\{f\circ \sigma\,\, ; \,\,
f\in{\mathcal
 H}({\mathscr S}_\alpha)\right\},
\]
with norm $\|f\circ \sigma\|_{\mathcal
H}=\|f\|_{{\mathcal
   H}({\mathscr S}_\alpha)}$, as follows from the equalities
\[
f(\sigma(\w))=\langle f(\cdot), K_{{\mathscr
S}_\alpha}(\cdot, \sigma(\w) )\rangle_{
{\mathcal H}({\mathscr S}_\alpha)}=
\langle f\circ \sigma (\cdot), K_{{\mathscr
S}_\alpha}(\sigma(\cdot),
\sigma(\w))\rangle_{\mathcal H}
\]
Thus the restrictions of the elements of
${\mathcal
 H}_2^\alpha({\mathbb D})$ to $\Omega_+$ are of the form as in
\eqref{hardy1} for $z$ restricted to $\Omega_+$.
By analytic  extension, the elements of
${\mathcal H}_2^\alpha({\mathbb D})$ have the
same
 form in the whole of ${\mathbb D}$.
\mbox{}\qed\mbox{}\\

\section{Schur multipliers}
\setcounter{equation}{0} A Schur function is a
function analytic and contractive in the open
unit disk. Equivalently, it is a function $s$
such that the operator of multiplication by $s$
is a contraction from the classical Hardy of the
open unit disk into itself. This last definition
is  our starting point to define
character-automorphic Schur multipliers.

\begin{Dn} A character-automorphic function $s^{\beta}(z)$, with character $\beta$, will be called a Schur
multiplier if the operator of multiplication by
$s^{\beta}(z)$ is a contraction from ${\mathcal
H}_2^{\overline{\beta}\alpha}({\mathbb D})$ into
${\mathcal H}_2^{\alpha}({\mathbb D})$.
\end{Dn}
Equivalently, the character-automorphic function
$s^{\beta}(z)$ is a Schur multiplier if and only
if the kernel
\begin{equation}
\label{multiplier} K^\alpha_{s^\beta}(z,w)=
k^\alpha(z,w)-s^{\beta}(z)s^{\beta}(w)^*k^{\overline{\beta}\alpha}(z,w)
\end{equation}
is positive in ${\mathbb D}$. The kernel
$K_{s^\beta}(z,w)$ is in particular positive in
$\Omega_+$, and we will consider it in
$\Omega_+$. We note that in view of \cite[Lemma
2, p. 142]{donoghue}, \cite[Theorem
 1.1.4, p. 10]{adrs}, the positivity of the analytic kernel $K_{s^\beta}$ on
$\Omega_+$ implies its positivity on $\D$.

\begin{Tm}
\label{schur} A character-automorphic function
$s^\beta$ is a Schur multiplier if and only if
there exists a ${\mathbb C}^{2\times 2}$--matrix
valued Schur function $\Sigma(z)$ such that
\begin{align}
s^\beta(z)&=\frac{A^\alpha(z)}{A^{\overline{\beta}\alpha}(z)}
\frac{\Sigma_{12}(\sigma(z))}{1-S_{{\overline{\beta}\alpha}}(z)\Sigma_{22}
(\sigma(z))}\\
S_\alpha(z)&=\Sigma_{11}(\sigma(z))+\frac{\Sigma_{12}(\sigma(z))
S_{\overline{\beta}\alpha}(z)\Sigma_{21}(\sigma(z))}
{1- S_{
\overline{\beta}\alpha}(z)\Sigma_{22}(\sigma(z))}
\end{align}
\end{Tm}
{\bf Proof:} The positivity of the kernel
\eqref{multiplier} in $\Omega_+$ is equivalent
to the positivity in $\Omega_+$ of the kernel
\[
K_{{\mathscr
S}_\alpha}(\sigma(z),\sigma(\w))-T(z)T(\w)^*
K_{{\mathscr
S}_{\overline{\beta}\alpha}}(\sigma(z),\sigma(\w)),
\]
where
\[
T(z)=\frac{A^{\overline{\beta}\alpha}(z)}{A^\alpha(z)}s^\beta(z).
\]
As in the proof of Proposition \ref{S0}, we note
that the function $\sigma$ maps $\Omega_+ \cap
\cF$ into some subset $\Delta \subset \D$ and
this mapping is one-to-one, and consider again
the function $\varsigma$ defined by: $(\varsigma
\circ \sigma) (z) = z, \, \forall z \in \Omega_+
\cap \cF$. The kernel
\begin{equation}
\label{newkernel} K_{{\mathscr
S}_\alpha}(\lambda,\mu)-T(\varsigma(\lambda))T(\varsigma(\mu))^*
K_{{\mathscr
S}_{\overline{\beta}\alpha}}(\lambda,\mu),
\end{equation}
is positive in $\Delta$. We now show that
$T\circ\varsigma$ admits an analytic extension
to ${\mathbb D}$. To that purpose we consider the
linear relation in ${\mathcal H}({\mathscr
S}_\alpha)\times {\mathcal H}({\mathscr
S}_{\overline{\beta}\alpha})$ spanned by the
elements of the form
\[
(K_{{\mathscr S}_\alpha}(\cdot,\w),
T(\varsigma(\w))^* K_{{\mathscr
 S}_{\overline{\beta}\alpha}}
(\cdot, \w)),\quad \w\in\Delta.
\]
It is densely defined and it is contractive
because of the positivity of the kernel
\eqref{newkernel} in $\Delta$. It is therefore
the graph of a densely defined contraction
$\widetilde{X}$. We note  its extension to
${\mathcal H}({\mathscr S}_\alpha)$ by $X$. For
$\w\in\Delta$ and $f\in {\mathcal H}({\mathscr
S}_\alpha)$ we have
\[
\begin{split}
(X^*f)(\w)&=\langle X^*f, K_{{\mathscr
S}_\alpha}(\cdot, \w)
\rangle_{{\mathcal H}({\mathscr S}_\alpha)}\\
&= \langle f, T(\varsigma(\w))^* K_{{{\mathscr
S}_{\overline{\beta}\alpha}}}(\cdot, \w)
\rangle_{{\mathcal H}({\mathscr S}_\alpha)}\\
&=T(\varsigma(\w))f(\w).
\end{split}
\]
Let $f_0(\lambda)=K_{{\mathscr
S}_{\overline{\beta}\alpha}}(\lambda, \w_0)$
where $\w_0\in{\mathbb D}$. We have
\[
T(\varsigma(\lambda))
=\frac{(X^*f_0)(\lambda)}{f_0(\lambda)},\quad
\lambda\in\Delta.\] It follows that
$T\circ\varsigma$ has an analytic extension to
${\mathbb D}$, which we will denote by
${\mathscr R}$. Thus the kernel
\begin{equation}
\label{newkernel1} K_{{\mathscr
S}_\alpha}(\lambda,\mu)-{\mathscr R} (\lambda)
{\mathscr R}(\mu)^* K_{{\mathscr
S}_{\overline{\beta}\alpha}}(\lambda,\mu),
\end{equation}
is analytic in $\lambda$ and $\mu^*$ in ${\mathbb
D}$. Therefore it is still positive in ${\mathbb
D}$; see \cite[Lemma 2, p. 142]{donoghue},
\cite[Theorem 1.1.4, p. 10]{adrs}. By
\cite[Theorem 11.1, p. 61]{ab6}, a necessary and
sufficient condition for the kernel
\eqref{newkernel1} to be positive is that there
exists a ${\mathbb C}^{2\times 2}$--matrix
valued Schur function $\Sigma(\lambda)$ such that
\begin{align}
\label{formulas} {\mathscr R}(\lambda)&=
\frac{\Sigma_{12}(\lambda)}{1-{\mathscr
S}_{\overline{\beta}\alpha}
(\lambda)\Sigma_{22}(\lambda)}\\
{\mathscr S}_\alpha(\lambda)&=
\Sigma_{11}(\lambda)+\frac{\Sigma_{12}(\lambda)
{\mathscr S}_{\overline{\beta}\alpha}
(\lambda)\Sigma_{21}(\lambda)} {1-{\mathscr
S}_{\overline{\beta}\alpha}
(\lambda)\Sigma_{22}(\lambda)}.
\end{align}
The above mentioned result from \cite{ab6} stems
from rewriting the positive kernel
\eqref{newkernel1} as
\[
\frac{{\mathcal A}(\lambda){\mathcal A}(\mu)^*-
{\mathcal B}(\lambda){\mathcal
B}(\mu)^*}{1-\lambda\mu^*}
\]
where
\begin{align*}
{\mathcal A}(\lambda)&=\begin{pmatrix}1&{\mathscr
R}(\lambda){\mathscr S}_{\overline{\beta}\alpha}(\lambda)\end{pmatrix}\\
{\mathcal B}(\lambda)&=\begin{pmatrix} {\mathscr
S}_\alpha(\lambda) &{\mathscr R}(\lambda)
\end{pmatrix},
\end{align*}
and using a factorization result known as
Leech's theorem, which insures the existence of
a ${\mathbb C}^{2\times 2}$-valued Schur function
$\Sigma$ such that
\[
{\mathcal B}(z)={\mathcal A}(z)\Sigma(z),
\]
from which \eqref{formulas} follows.\\

This unpublished result of R.B. Leech has been
proved using the commutant lifting theorem by M.
Rosenblum; see \cite[Theorem 2, p.
134]{rosenblum} and \cite[Example 1, p.
107]{rr-univ}. Further discussions and
applications can also be found in \cite{add}. It
can also be proved using tangential
Nevanlinna-Pick interpolation and
Montel's theorem.\\

Finally we replace in \eqref{formulas} $\lambda$
by $\sigma(z)$ where $z\in\Omega_+\cap{\mathcal
F}$. We obtain the formulas in the statement of
the theorem for $z\in \Omega_+\cap{\mathcal F}$,
and hence for $z\in{\mathbb D}$ by analytic
extension.
\mbox{}\qed\mbox{}\\

{\bf Acknowledgments.} The first author thanks
the Earl Katz family for endowing the chair
which supports his research. This work was done while the second author was visiting the department of mathematics of Ben-Gurion University. The second author thanks the Center for Advanced Studies in
Mathematics (CASM) of the department of
mathematics of Ben-Gurion University which supported his stay.
\bibliographystyle{plain}
\def\cprime{$'$} \def\lfhook#1{\setbox0=\hbox{#1}{\ooalign{\hidewidth
 \lower1.5ex\hbox{'}\hidewidth\crcr\unhbox0}}} \def\cprime{$'$}
 \def\cprime{$'$} \def\cprime{$'$} \def\cprime{$'$} \def\cprime{$'$}

\begin{minipage}[c]{\linewidth}
Daniel Alpay\\
Department of mathematics\\
Ben-Gurion University of the Negev\\
P.O.B. 653 Beer-Sheva 84105 - Israel\\
\textsf{E-mail: dany@cs.bgu.ac.il}
 \end{minipage}

\vspace*{1cm}

\begin{minipage}[c]{\linewidth}
Mamadou Mboup\\
UFR math\'ematiques et informatique \& EPI ALIEN INRIA\\
Universit\'e Paris Descartes\\
45, rue des Saints-P\`eres, 75006 Paris, France\\
\textsf{E-mail: mboup@math-info.univ-paris5.fr}
\end{minipage}
\end{document}